\documentclass[12pt]{article}
\usepackage{amssymb}
\usepackage{german}
\include{mak}
\newcommand{\du}[2]{\mbox{\raisebox{0.5ex}{$#1$}} 
\mbox{\hspace*{-0.25em}} \diagup_{\displaystyle \mbox{\hspace*{-0.3em}} #2}}
\begin{document}
\title{Representations of fundamental groups of algebraic manifolds 
and their restrictions to fibers of a fibration}
\author{J\"urgen Jost and Kang Zuo}
\maketitle

\section{Introduction}

We consider a surjective morphism $\, f: X\to Y\,$ from a smooth projective variety $X$ onto
a smooth projective variety $Y$ with connected fibers, henceforth called a fibration for short.
Typically, if a certain geometric object on $X$ like a cohomology class has a certain property,
then its restriction to a smooth fiber of $f$ trivially satisfies the same property. 
The converse question is of more interest: if the restriction of such a geometric object to a
smooth fiber enjoys a certain property, is this property also valid for the object on $X$ itself?
A prototype is the geometric version of the $(p,q)$ component theorem of Griffiths saying 
that if a class $H^k (X, {\mathbb C} )$ is of pure Hodge type $(p,q)$ at some smooth fiber 
$f^{-1} (y),$ then it has Hodge type $(p,q)$ at any smooth fiber.

In the so-called nonabelian cohomology, instead of classes in $H^1 (X, {\mathbb C} ),$ one considers
representations $\rho : \pi_1 (x) \to G$ into some linear algebraic group $G$. In the same way
as one associates a harmonic form to a cohomology class, one finds a $\rho$-equivariant harmonic
map $h:X \to \du{G}{K} $ into the symmetric space of noncompact type obtained as a homogeneous
space for $G$.  This harmonic map turns out to be pluriharmonic, meaning that its restriction 
to any subvariety is harmonic itself. We may thus reformulate the question indicated in the 
title of our paper, namely what one can infer about a representation of $\pi_1 (X)$ if one knows
a relevant property of the induced representation on $\pi_1 (f^{-1} (y) ),$ or, more generally,
of the one on $\pi_1 (Z), Z$ a generic subvariety of $X,$ as the question of what we can deduce
about a pluriharmonic map from its restriction to $f^{-1} (y)$ or $Z$. Let us start with some 
easy observations in this direction before formulating our actual results. A harmonic map into
$\du{G}{K}$ is constant on any rational variety. On one hand this implies that a pluriharmonic
map $h$ into $\du{G}{K}$ is an invariant of the birational class of $X$, as has been observed
by many people, but on the other hand this tells us that we cannot deduce information about 
$h$ from its restriction to rational subvarieties. In a more positive direction, if we have a 
smooth map $g:X \to \du{G}{K} $ whose restriction to any fiber $f^{-1} (y)$ is harmonic, and
if we are in a situation where the image of the fibers is sufficiently large so that harmonic
map uniqueness holds (this is the case if there are no parallel vector fields along the image),
then $g$ itself must be harmonic. In particular, in that situation if we have homomorphisms 
$\rho : \pi_1 (f^{-1} (y)) \to G$ on the fundamental groups of the fibers, then these 
homomorphisms induce a unique homomorphisms $\rho : \pi_1 (X) \to G.$ 
An essential point of the present note consists in results that refine
the preceding simple observation. We shall only need the weaker assumption that the harmonic
maps on the fibers are nonconstant, or equivalently that the representations on 
$\pi_1 (f^{-1} (y))$ are nontrivial, plus the assumption that the image of $\pi_1 (X)$ is 
Zariski dense. Then the representation of $\pi_1 (X)$ cannot be deformed with the representation
on $\pi_1 (f^{-1} (y))$ kept fixed (see Thm. 2b below).
We also study a similar situation where instead of the fibers of a holomorphic map we consider
a subvariety $Z$ whose fundamental group surjects onto the one of $X.$ Again, in that case, the 
representation on $Z$ locally determines the one on $X$ (see Thm. 2a). 
Finally, we study how other properties like coming from a variation of Hodge structures or 
arithmeticity on the fibers of a holomorphic map induce the corresponding ones on $X$ 
(Thm. 1 and 2c, resp.).\\

{\bf Remark} We wrote up the first version of this note in December 1999. After we had circulated 
our result, also a preprint by Pantev and Katzarkov appeared where they proved somewhat weaker results
by a different method. This result was announced by them some time ago.\\

{\bf Acknowledgments} During the preparation of this paper, the second author was supported by a Heisenberg
fellowship of the DFG.

\section{Results and proofs}
%

{\bf Theorem 1}  {\sl  Let $\, f: X\to Y\,$ be a fibration and 
$\, Z=\cup_{i=1}^mZ_i\,$ be a reduced fibre of $\, f\,$ with smooth irreducible components.
 Suppose $\, \rho:\pi_1(X)\to G\,$ is a Zariski dense representation into an almost 
simple algebraic group $\, G\,$ and does not factor through $\, f.\,$ 
If the restriction $\, \rho|_{Z_i},\,1\leq i\leq m\,\,$ comes from $\, {\mathbb Z}-$variations 
of Hodge structure, then $\, \rho\,$ itself comes from $\, {\mathbb Z}-$variations of Hodge structure.}

We say a representation $\, \rho:\, \pi_1(X)\to G\,$ factors through a fibration $\,f:X\to Y,\,$ 
if there exists a finite etale covering with a blowing up $\,e: X'\to X\,$ and a representation 
$\, \tau:\, \pi_1(Y')\to G,\,$ where $\, Y'\,$ is the base of the Stein-factorisation 
$\, f':X'\to Y'\,$ of $\,fe: X'\to X\to Y,\,$ such that $\, e^*(\rho)=f'^*(\tau).\,$ 

{\bf Remark 1} The following consideration was pointed out to the authors by C. Simpson.
 We can  replace $\, \rho\,$ by a section of the  relative Betti space $\,
M_B(X/Y, G)\,$ which is flat with respect to  the nonabelian Gau\3-Manin connection
introduced by Simpson [S2]. We shall prove Theorem 1 for such a flat section in a forthcoming 
version.\\
A flat section glues to a  global representation $\, \tilde \rho'\,$ 
on the covering   $\,\tilde X'\to X\,$ corresponding to 
 $\,\pi_1(X)\to \pi_1(Y)\to 1.\,$ The decktransformation group $\, \pi_1(Y)\,$ acts 
on $\, M_B(\tilde X',G).\,$ One shall try to show that $\, \tilde \rho'\,$ is
$\,\pi_1(Y)-$invariant and  descends to a representation $\, \rho:\pi_1(X)\to G.\,$

Theorem 1 is an easy consequence of Simpson's theorem about when a representation will
come from variations of Hodge structure and Theorem 2 below, which can be considered as a 
type of  Lefschetz hyperplane theorem for possibly singular subvarieties, whose fundamental 
groups surjects onto $\,\pi_1(X),\,$  or for subvarieties arising as fibres of fibrations. 

{\bf Theorem} (Simpson, Corollary 4.2, [S1]) {\sl The representations of $\, \pi_1(X)\,$ 
which come from complex variations of Hodge structure are exactly the semisimple ones 
which are fixed by the action of $\, {\mathbb C}^*$.}

Hitchin originally defined this action in the form of an action of $\, U(1)\subset 
{\mathbb C}^*\,$ [H].

Let $\, Z=\cup_{i=1}^m Z_i\subset X\,$ be a subvariety and suppose that all irreducible
components are smooth. Given a representation $\, \rho\in M_B(X,G),\,$ the pull back via 
$\, Z_i\hookrightarrow X\,$ defines a point $\, \rho|_{Z_i}\in M_B(Z_i,G).\,$ 
Hence, that defines a morphism

               $$ r:  M_B(X, G)\to \prod_{i=1}^m M_B(Z_i, G).$$

We may consider $\, G\subset SL_n\,$ as a group scheme defined over some number field and with 
a fixed integral structure (for instance, induced by $\, SL_n\,$).  So, the morphism $\,r\,$  
is clearly also defined over some number field. 

{\bf Theorem 2} {\sl {\bf a)} Suppose that  the homomorphism $\, \pi_1(Z)\to \pi_1(X)\,$ is surjective. 
Then the preimage of $\, r\,$ over any point  consists of finitely many points only.}

{\sl {\bf b)} Suppose that $\, Z=\cup Z_i\,$ is a reduced fibre of a surjective morphism 
$\, f: X\to Y\,$ with connected fibres, and that $\, \rho:\,\pi_1(X)\to G\,$ is Zariski dense 
and does not factor through $\, f.\,$ Then $\,\rho\in r^{-1}( r(\rho))\,$ is an isolated point.}

{\sl {\bf c)} Let $\, Z=\cup_{i=1}^m Z_i\,$ and $\, \rho\,$ be the same as in b). If the restriction 
$\, \rho|_{Z_i},\, 1\leq i\leq m\,$ is valued in the ring of the algebraic integers of a number field, 
then this also holds for $\,\rho\,$ itself.}

{\sl Proof of Theorem 1} \quad Applying Simpson's theorem, the action of 
$\, {\mathbb C}^*$ fixes $\, r(\rho).\,$ Since the action of 
$\, {\mathbb C}^*$ fixing $\, r(\rho))\,$ commutes with $\, r,\,$ the fibre 
$\, r^{-1}(r(\rho)\,$ is fixed by the $\, {\mathbb C}^*-$action. Since 
$\, \rho\in  r^{-1}(r(\rho))\,$ is an isolated point by b) in Theorem 2, $\, \rho\,$ is a fix point of
the $\, {\mathbb C}^*-$action. Applying Simpson's theorem again, $\, \rho\,$ comes from complex 
variations of Hodge structure. The integral property of $\, \rho\,$ follows from c) in Theorem 2.

{\sl Proof of Theorem 2 a)}\quad  Since the morphism 

$$ r:  M_B(X, G)\to \prod_i M_B(Z_i, G) $$   
 
is defined over some number field, we only need to check the property a) for any point 
$\,\tau=(\tau_1,...,\tau_m)\in \prod_{i=1}^m M_B(Z_i, G),\,$ valued in some number field $\, K.\,$

Since $\, \pi_1(Z_i), 1\leq i\leq m\,$ is a finitely presented group, we may find a prime ideal $\, p\,$ 
of $\, {\cal O}_K\,$ such that $\, \tau_i(\pi_1(Z_i))\subset G( {\cal O}_{K_p}),1\leq i\leq m,\,$ where 
$\,K_p\,$ is the local field at the place $\, p.\,$  

Consider now the morphism $\, r\,$ defined over $\, K_p.\,$ If the statement a) were not true, then 
$\, r^{-1}(\tau)\,$  would contain a positive dimensional component $\, C.\,$ Hence, 
we may find some representation $\,\rho\in C\,$ which is valued in some finite extension of $\, K_p\,$ 
and is $p-$unbounded. Notice that $\, r(\rho)\,$ is $p-$bounded.

Now let 

               $$ u_{\rho}:\tilde X\to \triangle (G(L_p))$$ 

denote the $\rho-$equivariant pluriharmonic map into the corresponding Bruhat-Tits building. 
The existence of such harmonic maps is shown in [GS]. Since $\, \rho\,$ is $p-$unbounded, 
$\, u_{\rho}\,$ is not constant. 

On the other hand, the surjectivity of $\, \pi_1(Z)\to \pi_1(X)\,$ and the $p-$boundedness of 
$\, r(\rho)\,$ imply that $\, u_{\rho}\,$ is constant. The following argument can be found in [LR].

Since $\, \pi_1(X)\to \pi_1(Z)\,$ is surjective, the preimage $\, \tilde Z=\cup \tilde Z_i\,$
is connected.  The restriction $\, u_{\rho}|_{\tilde Z_i}\,$ is just the corresponding equivariant 
harmonic map of $\, \tau_i\,$ and is constant, since $\, \tau_i\,$ is $p-$bounded. 
Hence, $\, u_{\rho}(\tilde Z)\,$ is a point $\,q.\,$ Therefore, the action of $\, \rho(\pi_1(X) )\,$ 
also fixes $\,q.\,$ This implies that $\, u_{\rho}(\tilde X)=q.\,$ A contradiction.

{\sl Proof of  b)}\quad  We may find  a stratification on 

$$ r: M_B(X,G)\to \prod_{i=1}^mM_B(F_i,G)$$ 

which is defined over some number field and such that:

i) $\, r(M_B(X,G)_j)= (\prod_{i=1}^mM_B(F_i,G))_j\,$ and  

ii) $\, r:M_B(X,G)_j\to (\prod_{i=1}^mM_B(F_i,G))_j\,$ is flat.

We want to show that if $\,\rho\in r^{-1}( r(\rho))\,$ is not  an isolated point, then $\, \rho\,$
factors through $\, f\,$ after passing to a finite etale covering and a blowing up $\, X'\to X.\,$

Using the above stratification we may first  show that property  for  those $\, \rho',\,$  who are 
in the same strata as $\,\rho\,$ and valued in some number field. If all such $\, \rho'\,$ factor 
through $\, f,\,$ then $\, \rho\,$ also factors through $\,f.\,$

Suppose $\, \rho\,$ is valued in some number field $\, K.\,$ By the same reason as explained in the 
proof of a), we may find a local field $\, K_p,\,$ such that $\, r(\rho)\,$ is valued in 
$\, {\cal O}_{K_p}.\,$

If $\, \rho\in r^{-1}(r(\rho))\,$ were not an isolated point, then we would find an irreducible curve 
$\, \rho \in C\subset  r^{-1}(r(\rho))\,$ that contains infinitely many $p-$unbounded and Zariski dense 
representations $\,\rho_i: \pi_1(X)\to  G(L_{p,i}),\,$ where $\, L_{p,i}\,$ is a finite extension of 
$\, K_p.\,$

Let

$$ u_{\rho_i}:\tilde X\to \triangle (G(L_{p,i}))$$ 

denote the $\rho_i-$equivariant pluriharmonic map into the Bruhat-Tits building.  
$\,u_{\rho_i}\,$ is not constant, since $\, \rho_i\,$ is not $p-$bounded. 

Consider the pulled back fibration 

$$ \tilde f: \tilde X\to \tilde Y.$$

We have  

{\bf Claim 1} {\sl $\, u_{\rho_i}\,$ factors through  $\,\tilde f.\,$ }

{\sl Proof of Claim 1} \quad  The differential $\, d'u_{\rho_i}\,$ is a collection of holomorphic 
1-forms $\, \theta_i\,$ on a finite ramified covering $\, X^s\to X.\,$
The pull back of $\, \theta_i\,$ to the corresponding fibre $\, F^s=\cup_{i=1}^m F^s_i\,$ is zero, 
since $\, r(\rho)\,$ is $p-$bounded and $\, u_{\rho_i}|_{F_i},\,1\leq i\leq m\,$ is constant. 

The pull back of $\, \theta_i\,$ to any fibre of $\, f^s\,$ is also zero, since any closed 1-cycle on a 
fibre of $\, f^s\,$ is homotopic to some 1-cycle in $\, F^s\,$ and  the integration of $\, \theta_i\,$ 
on any closed 1-cycle on a fibre of $\, f^s\,$ is zero. 
Therefore, we see that $\, u_{\rho_i}\,$ factors through $\,\tilde f.\,$

{\bf Claim 2} {\sl  Let  $\,F_0:= f^{-1}(y_0)\,$ be a smooth fibre of $\, f.\,$ Then $\, 
\rho_i(\pi_1(F_0))\,$ is a finite subgroup of $\, G(K_p).\,$}

{\sl  Proof of Claim 2} \quad  Pulling back $\, f\,$ to the universal coverings, by Claim 1 the 
harmonic map $\, u_{\rho_i}\,$ factors through $\, \tilde f,\,$

\[ \begin{array}{ccc}
\tilde X &\stackrel{\tilde f}{\vector(1,0){30}}&\tilde Y\\
&  {\vector(1,-1){30}}&{\vector(0,-1){30}} \\
&     & \triangle(G(K_p))\end{array}\]

Fixing a base point $\, x_0\in F_0,\,$ let  $\, \Gamma=:im(\pi_1(F_0,x_0)\to \pi_1(X,x_0)).\,$ 
We first want to show that $\, \Gamma\,$ fixes a unbounded subset in $\, \triangle(G(L_{p,i})).\,$\\
Let $\,  X_0\subset X\,$ be the Zariski open set, such that the map $\, f:\, X_0\to Y_0\,$ is regular. 
We denote by $\, \tilde X_0\to X_0\,$ the universal covering of $\,\tilde X_0.\,$

The subgroup $\, \Gamma\subset \pi_1(X_0,x_0) \,$ operates on a connected component $\, \tilde F_{0,0}\,$ 
of the preimage $\, \tilde F_0\subset \tilde X_0.\,$ Since the harmonic map $\, u_i\,$ is 
$\rho-$equivariant and factors through the fibration $\, \tilde f:\tilde X_0\to \tilde Y_0\,$ to an 
equivariant harmonic map

$$ v_i:\tilde Y_0\to \triangle(G(L_{p,i})),$$

$\,\rho_i(\Gamma)\,$ fixes the image $\, u_i(\tilde F_{0,0})=v_i(\tilde y_{0,0})=: z_{0,0}.\,$\\
If $\, \tilde F_{0,j}\,$ is another connected component of $\, \tilde F_0,\,$
then there exists an element $\, g_j\in\pi_1(X_0,x_0),\,$ such that $\, g_i(\tilde F_{0,0})=\tilde 
F_{0,j}.\,$ So, the conjugation $\, g_j\Gamma g^{-1}_j\,$ operates on $\, \tilde F_{0,j},\,$ 
and by the same reason as above $\,\rho_i(g\Gamma g^{-1})\,$ fixes $\, u_i(\tilde F_{0,j})
=v_i(\tilde y_{0,j}):=z_{0,j}.\,$\\
Considering the exact sequence of the homotopy groups (coming from the definition of $\Gamma$)

$$ 1\to\Gamma\to \pi_1(X_0,x_0)\to \pi_1(Y_0,y_0)\to 1,$$

we see in particular that $\, \Gamma\subset \pi_1(X_0,x_0)\,$ is a normal subgroup. 
Hence, $\,\rho_i (\Gamma)= \rho_i(g_j\Gamma g^{-1}_j)\,$ fixes $\, z_{0,j},\forall j\in J.\,$\\

Now we show that the subset $\, Z:=\{z_{0,j}\}_{j\in J}\subset\triangle(G(L_{p,i})) \,$ is unbounded.
Since each  different point in the subset $\, \{y_{0,j}\}_{j\in J}\subset \tilde Y_0\,$ is contained in 
a different fundamental domain $\, D_j\subset \tilde Y_0\,$  and the images
$\, v_i(D_j)\ni z_{0,j}, \,j\in J\,$ are uniformly bounded (they are permuted by 
$\, \rho_i(\pi_1(X_0,x_0))\,$ as isometry), together with the unboundedness of $\, v_i(\tilde Y_0)\,$ 
this implies that $\, Z\,$ must be an unbounded subset in $\,\triangle(G(L_{p,i})).\,$

We want to show further that $\, \rho_i(\Gamma)\,$ fixes a point on the boundary of 
$\,\triangle(G(L_{p,i})).\,$\\
Since $\, Z\subset \triangle(G(L_{p,i}))\,$ is an unbounded subset, the convex subcomplex generated by 
$\, Z\,\,$ contains at least one geodesic line $\, L.\,$ Since $\, \rho_i(\Gamma)\,$ fixes $\, Z\,$ 
pointwisely, $\, \rho_i(\Gamma)\,$ fixes this convex subcomplex pointwisely. Hence, 
$\,\rho_i( \Gamma)\,$ fixes $\, L\,$ and the point on the boundary of $\, \triangle(G(L_{p,i})\, $ 
defined by $\, L.\,$ \\
That shows that $\, \rho_i(\Gamma)\,$  is contained in a parabolic subgroup of $\, G.\,$ 
In particular, it is not Zariski dense in $\,G(L_{p,i}).\,$ Furthermore, the exact sequence of the 
homotopy groups  above shows that the Zariski closure of $\, \rho_i(\Gamma)\,$ is a
normal  algebraic subgroup in $\, G(L_{p,i}).\,$ Since $\, G(L_{p,i})\,$ is almost simple,
$\,\rho_i( \Gamma)\,$  must be finite.  Claim 2 is proved.


By restriction to $\, F_0\,$ we get a family of representations 

$$ \rho_t :\pi_1(F_0)\to G(K_p),\quad t\in C.$$

By Claim 2 $\, \rho_i(\pi_1(F_0))\subset G(K_p),\, i\in I\,$ is not Zariski dense.
Since $\, I\,$ is infinite and $\, C\,$ is irreducible, the subset $\, \{\rho_i|_{F_0}\}_{i\in I}\subset C\,$
is Zariski dense. Since $\, G\,$ is an almost
simple group,  the Zariski density of representations is a Zariski open condition in the  space of representations. That shows
$\, \rho(\pi_1(F_0))\subset G(K_p)\,$ is not Zariski dense, too. By the same reason as in Claim 2, 
$\, \rho(\pi_1(F_0))\subset G(K_p)\,$ is a finite subgroup.

Since $\, \rho (\pi_1(X)\subset G\,$ is residually finite, one may find a finite etale covering with a 
blowing up $\, X'\to X\,$ such that the pull back of $\, \rho\,$ factors through the Stein 
factorisation of $\, f:X'\to Y.\,$ So, this leads to a contradiction to our assumption in b).

{\sl Proof of c)}  Since $\, r(\rho)\,$ is valued in some number field $\, K\,$ and $\, \rho\in r^{-1}( r(\rho))\,$
is an isolated point, $\,\rho\,$ is valued in some finite extension  $\,L\supset K.\,$
Now as $\,r(\rho)\,$ is bounded with respect to any prime ideal of $\, {\cal O}_K,\,$ the same argument
in b) shows that $\, \rho\,$ is also  bounded with respect to any prime ideal of $\, {\cal O}_L.\,$ 


\section{References}

\lbrack G \rbrack \quad Griffiths, P.: Periods of integrals on algebraic manifolds III. 
{\sl Publ. Math. IHES 38 (1970), 125-180.}\\

\lbrack GS\rbrack \quad Gromov, M. and Schoen, R.: Harmonic maps into singular spaces and $p-$adic 
superrigidity for lattices in groups of rank one.  {\sl Publ. Math. IHES 76 (1992), 165-246.}\\

\lbrack H\rbrack \quad Hitchin, N.J.: The self-duality equations on a Riemann surface, {\sl Proc.
 London Math. Soc. (3), 55 (1987), 59-126.}\\

\lbrack LR\rbrack \quad Lasell, B. and Ramachandran, M.: Observations on harmonic maps and singular
varieties. {\sl Ann Sci. Ecole Norm. Sup. 29. (1996), 135-148.}\\

\lbrack S1\rbrack \quad Simpson, C.: Higgs bundles and local systems.  
{\sl Publ. Math. IHES 75 (1992), 5-95.}\\

\lbrack S2\rbrack \quad Simpson, C.: The Hodge filtration on nonabelian cohomology. 
{\sl Proc. Symp. Pure Math., Vol 62, Part 2, 217-281.}

\end{document}